\documentclass[letterpaper, 11pt]{amsart}

\usepackage[all]{xy}        % Commutative Diagrams
\CompileMatrices            % Faster Tex-ing
\UseTips                    % Use Computer Modern Arrowheads

\usepackage[german,english]{babel}   %Foreign language support
\usepackage[bookmarks=true]{hyperref}   % Hyperref in DVI and PDF (like HTML Links between sections)
\usepackage{amssymb}
\usepackage{xspace}
\theoremstyle{plain}
\newtheorem{theorem}{Theorem}[section]
\newtheorem*{theorem*}{Theorem}

\newtheorem{lemma}[theorem]{Lemma}

\theoremstyle{definition}

\theoremstyle{remark}
\newtheorem{remark}[theorem]{Remark}
\newtheorem{example}[theorem]{Example}

\numberwithin{equation}{section}

\newcommand{\enm}[1]{\ensuremath{#1}}          % Shortcuts
\newcommand{\op}[1]{\operatorname{#1}}
\newcommand{\cal}[1]{\mathcal{#1}}

\newcommand{\CC}{\enm{\mathbb{C}}}             % All Number domains easily accssable

\newcommand{\QQ}{\enm{\mathbb{Q}}}
\newcommand{\ZZ}{\enm{\mathbb{Z}}}
\newcommand{\FF}{\enm{\mathbb{F}}}

\newcommand{\PP}{\enm{\mathbb{P}}}

\newcommand{\Dd}{\enm{\cal{D}}}

\newcommand{\Oo}{\enm{\cal{O}}}

\renewcommand{\phi}{\varphi}        % Dont know how to not loose the original ones???
\renewcommand{\theta}{\vartheta}
\renewcommand{\epsilon}{\varepsilon}

\newcommand{\Spec}{\op{Spec}}
\newcommand{\Proj}{\op{Proj}}

\newcommand{\Ext}{\op{Ext}}

\newcommand{\Supp}{\op{Supp}}

\newcommand{\tensor}{\otimes}         % Symbols with meaning

\renewcommand{\to}[1][]{\xrightarrow{\ #1\ }}

\newcommand{\usc}[1][m]{\underline{\phantom{#1}}}

\newcommand{\defeq}{\stackrel{\scriptscriptstyle \op{def}}{=}}

\newcommand{\Sol}{\op{Sol}}

\newcommand{\Rb}{\mathbf{R}}

\newcommand{\et}{\text{\'et}}

\newcommand{\draft}[1]{}

\begin{document}

%
%   Preliminary information     %%%%%%%%%%%%%%%%%%%%%%%%%%%%%%%%%%%%%%%%
%

\title[]{Lyubeznik's invariants for cohomologically isolated singularities.}
\author{Manuel Blickle}
\address{Universit\"at Essen, FB6 Mathematik, 45117 Essen,
Germany} \email{manuel.blickle@uni-essen.de} \urladdr{\url{www.mabli.org}}
%\author{Raphael Bondu}
%\address{Universit\"at Essen, FB6 Mathematik, 45117 Essen,
%Germany} \email{raphael.bondu@uni-essen.de}
\thanks{The author is supported by the DFG Schwerpunkt \emph{Globale
Methoden in der komplexen Geometrie}.}
%
%   Title and abstract (table of contents)     %%%%%%%%%%%%%%%%%%%%%%%%%
%

\maketitle

\begin{abstract}
    In this note I give a description of Lyubeznik's local cohomology
    invariants for a certain natural class of local rings, namely the ones which have
    the same local cohomology vanishing as one expects from an isolated singularity.
    This strengthens our results of \cite{BliBon.LocCohomMult} while at the
    same time somewhat simplifying the proofs. Through examples I further point out the
    bad behavior of these invariants under reduction mod $p$.
\end{abstract}

\section{Introduction}
Let $A=R/I$ for $I$ an ideal in a regular (local) ring $(R,m)$ of dimension $n$
and containing a field $k$. The main results of
\cite{Lyub.FinChar0,HuSha.LocCohom} state that the local cohomology module
$H^a_m(H^{n-i}_I(R))$ is injective and supported at $m$. Therefore it is a
finite direct sum of $e=e(H^a_m(H^{n-i}_I(R)))$ many copies of the injective
hull $E_{R/m}$ of the residue field of $R$. Lyubeznik shows in
\cite{Lyub.FinChar0} that this number
        \[
        \lambda_{a,i}(A) \defeq e(H^a_m(H^{n-i}_I(R)))
        \]
does not depend on the auxiliary choice of $R$ and $I$. If $A$ is a
complete intersection, these invariants are essentially trivial (all
are zero except $\lambda_{d,d}=1$ where $d=\dim A$). The goal of
this paper is to describe these invariants for a large class of
rings, including those which are complete intesections away from the
closed point. Alternatively this class of rings can be viewed as
consisting of the rings which behave cohomologically like an
isolated singularity.
\begin{theorem}\label{thm.main}
Let $A=\Oo_{Y,x}$ for $Y$ a closed $k$--subvariety of a smooth variety $X$. If
for $i \neq d$ the modules $H^{n-i}_{[Y]}(\Oo_X)$ are supported in the point
$x$ then
\begin{enumerate}
    \item For $2 \leq a \leq d$ one has
    \[
        \lambda_{a,d}(A)-\delta_{a,d} = \lambda_{0,d-a+1}(A)
    \]
    and all other $\lambda_{a,i}(A)$ vanish.\draft{\footnote{This vanishing condition is
    essentially equivalent to the assumption on the support.
    Clearly if $H^{n-i}_I(R)$ is supported in $m$ for $i \neq d$
    then for $a \neq 0$ and $i \neq d$, $\lambda_{a,i}=0$
    by definition of $\lambda_{a,i}$. The converse is slightly more involved.}}
    \item
    \[
    \lambda_{a,d}(A)-\delta_{a,d} =
    \begin{cases}
        \dim_{\FF_p} H^{d-a+1}_{\{x\}}(Y_{\et},\FF_p) & \text{if $\op{char} k = p$} \\
        \dim_{\CC}   H^{d-a+1}_{\{x\}}(Y_{\op{an}},\CC) & \text{if $k = \CC$}
    \end{cases}
    \]
\end{enumerate}
where $\delta_{a,d}$ is the Kroneker delta function.
\end{theorem}
In the case that $A$ has only an isolated singularity and $k=\CC$,
this was shown by Garcia~Lopez and Sabbah in \cite{LopezSabbah}. In
\cite{BliBon.LocCohomMult} Bondu and myself proved part (1) in all
characteristics but part (2) only under an additional assumption.
The point of this note is to show that this additional assumption
was unnecessary, an observation which also makes the proof much
clearer.

\section{Proof of the result}
Let us first fix our notation. Let $Y \subseteq X$ be a closed subscheme of
$X$. Let $X$ be smooth of dimension $n$ and $Y$ is of dimension $d$. Let $x \in
Y$ be a point. We take the freedom to shrink $X$ (and $Y$) since everything is
local at $x$. We denote the inclusions as follows:
\[
\xymatrix{ {Y} \ar[r]^{i'} & {X} \\
           {Y-x} \ar[u]_j\ar[r]^i & {X-x} \ar[u]_{j'} }
\]
To ease notation we carry out the argument in the case when $\op{char} k$ is
positive. The proof of the characteristic zero result is exactly the same,
replacing the Emerton--Kisin correspondence by the Riemann--Hilbert
correspondence.

Building upon \cite{BliBon.LocCohomMult} it remains to show that
\[
    \lambda_{a,d}(A)-\delta_{a,d}= \dim H^{d-a+1}_{\{x\}}(Y_{\et},\FF_p)
\]
for $2 \leq i \leq d$. By definition of the invariants
\[
    \lambda_{a,i}(A) = e(H^a_{[x]}(H^{n-i}_{[Y]}(\Oo_X)))
\]
this comes down to computing $\dim \Sol(H^a_{[x]}(H^{n-i}_{[Y]}(\Oo_X))$. The
only trick in the proof is to replace $H^{n-i}_{[Y]}(\Oo_X)$ by something more
accessible -- that is by something whose solutions $\Sol$ can readily be
computed. The rest is mere computation.

Quite generally one has a short exact sequence
\[
    0 \to K \to j'_{!*}H^{n-i}_{[Y]}(\Oo_X)|_{X-x} \to
    H^{n-i}_{[Y]}(\Oo_X) \to C \to 0
\]
where, by construction, the kernel $K$ and quotient $C$ are both supported on
$x$.\footnote{The functors we call $j_*$ and $j_{!*}$ are denoted by $j_+$ and
$j_!+$ in \cite{EmKis.Fcrys}. In order to adjust to the notation used in the
Riemann--Hilbert correpondence I changed this notation here.}\footnote{The
intermediate extension $j_{!*}M$ of a module ($\Dd_{(X-x)}$-- or
$\Oo_{F,(X-x)}$--module if in characteristic $0$ or $p$ respectively) is
defined as the unique smallest submodule $M'$ of $H^0j_*(M)$ such that
$j^!M'\cong M$, where $j: X-x \to X$ is the open inclusion of the complement of
$x$ into $X$. The intermediate extension appears here as a substitute for $j_!$
which does not exist in the characteristic $p$ context. It was this realization
(replace $j_!$ with $j_{!*}$) which made it possible to make the argument work
in all characteristics.} Splitting this four term sequence into two short exact
sequences and using the long exact sequence for $\Gamma_{[x]}(\usc)$ we get
that for $a \geq 2$
\begin{equation}\label{eq}
    H^a_{[x]}(H^{n-i}_{[Y]}(\Oo_X)) \cong
    H^a_{[x]}(j'_{!*}(H^{n-i}_{[Y]}(\Oo_X)|_{X-x})).
\end{equation}
This first substitution, combined with \cite[Lemma
2.3]{BliBon.LocCohomMult}, we record as a Lemma:
\begin{lemma}\label{lem.1}
    With notation as above, and without any assumptions on the
    singularities of $A=\Spec \Oo_{Y,x}$
    \[
        \lambda_{a,d}(A)= \dim (H^{-a}j'_{!*}\Sol
        H^{n-d}_{[Y-x]}(\Oo_{X-x}))_x
    \]
    for $2 \leq a \leq d$.
\end{lemma}
\begin{proof}
By definition we have
\[
\begin{split}
    \lambda_{a,d}(A)&=e(H^a_{[x]}H^{n-d}_{[Y]}(\Oo_X)) \\
                    &=\dim (\Sol H^a_{[x]}H^{n-d}_{[Y]}(\Oo_X)))_x \\
                    &=\dim (k_{!}k^{-1}H^{-a}\Sol
                    H^{n-d}_{[Y]}(\Oo_X)))_x\qquad
                    \text{\cite[Lemma 2.3]{BliBon.LocCohomMult}}
                    \\
                    &=\dim (H^{-a}\Sol j'_{!*}(H^{n-d}_{[Y]}(\Oo_X)|_{X-x}))_x
\end{split}
\]
where $k: x \to X$ is the inclusion of the point. The commutation of $\Sol$
with $j'_{!*}$ now finishes the argument.
\end{proof}

The assumption that $H^{n-i}_{[Y]}(\Oo_X)$ is supported at $x$ for $i \neq d$
we rephrase by saying that one has a quasi-isomorphism of complexes
\[
    H^{n-d}_{[Y-x]}(\Oo_{X-x}) \cong
    \Rb\Gamma_{[Y-x]}(\Oo_{X-x})[n-d].
\]
and the solutions of the latter can be computed easily\footnote{In the case
$k=\CC$ this computation yields $i_!\CC_{Y-x}[d]$ instead.} (and is preverse!)
as done in \cite{BliBon.LocCohomMult}:
\[
    \Sol(\Rb\Gamma_{[Y-x]}(\Oo_{X-x})[n-d]) = i'_!(\FF_p)_{Y-x}[d]
\]
Since $j$ is just the inclusion of complement of a point we have that
\[
    j_{!*}(\usc) \cong \tau_{\leq {d - 1}} \Rb j'_*(\usc)
\]
by \cite[V.2.2 (2)]{Borel.Icohom}. Continuing the computation of Lemma
\ref{lem.1} using these observations we get for $a \geq 1$ that
\[
\begin{split}
    \lambda_{a,i}(A) &= \dim (H^{-a}j'_{!*}\Sol H^{n-d}_{[Y-x]}(\Oo_{X-x}))_x \\
                     &= \dim (H^{-a}j'_{!*}\Sol (\Rb\Gamma_{[Y-x]}(\Oo_{X-x})[n-d]))_x \\
                     &= \dim (H^{-a}j'_{!*}i'_!(\FF_p)_{Y-x}[d])_x \\
                     &= \dim (H^{-a}i_!j_{!*}(\FF_p)_{Y-x}[d])_x \\
                     &= \dim (H^{d-a}\tau_{\leq d-1}\Rb j_*(\FF_p)_{Y-x})_x \\
                     &= \dim (H^{d-a}\Rb j_*(\FF_p)_{Y-x})_x
\end{split}
\]
The latter was computed in \cite[Lemma 2.7]{BliBon.LocCohomMult} to be equal to
$H^{d-a+1}_{\{x\}}(Y_\et,\FF_p)+\delta_{a,d}$ as required.
\begin{remark}
    A slight refinement of the same techniques yield also to a more general
    statement. Namely, if one requires vanishing of $H^a_{m}H^{n-i}_I(R)|_{\Spec R-{\op{pt}}}$ below
    the diagonal $a=i+m$ then the result remains true in the range
    $d-m+2 \leq a \leq d$.
\end{remark}

\section{Examples}
This section is to provide some examples of the bad behavior of the
invariants $\lambda_{a,i}$ under reduction to positive
characteristic. The uniformity of the Theorem for all
characteristics seems, on the first sight, to suggest that that one
can expect a good behavior of the invariants under reduction mod
$p$. This impression is however quickly shattered, essentially for
reasons that local cohomology is well known to behave poorly under
reduction. Alternatively one also can observe that the cohomology
theory corresponding to $H^i_{x}(Y_{an},\CC)$ under reduction is not
$H^i_{x}(Y_\et,\FF_p)$ but rather $p$-adic rigid cohomology or
crystalline cohomology, of which $H^i_{x}(Y_\et,\FF_p)$ is only a
small part, namely the slope zero part.

The examples that follow are standard examples for the bad behaviour
of local cohomology under reduction mod $p$. I learned them from a
talk by Anurag Singh at Oberwofach in March 2005. Our general setup
is as follows: Let $A = R/I$ where $R$ is a polynomial ring over
$\ZZ$ and $I$ is a homogeneous ideal. We denote by $A_0=A
\tensor_\ZZ \CC$ the generic characteristic zero model and by $A_p =
A \tensor \FF_p$ for all $p$ prime the positive characteristic
models.
\begin{example}
Let $R=\ZZ\left[\begin{matrix} u & v & w \\ x & y & z
\end{matrix}\right]$ and $I=(\delta_1,\delta_2,\delta_3)$ be the
ideal of $2 \times 2$ minors of the displayed matrix of variables.
Then $A = R/I$ has a free resolution (as an $R$--module) given by
\[
    0 \to R^2 \to[\scriptstyle\begin{pmatrix} u & x \\ v & y \\ w & z \end{pmatrix}] R^3 \to[(\delta_1\; \delta_2\; \delta_3)] R \to 0
\]
This shows that $\Ext^3(R/I,R)=0$ and therefore, reducing mod $p$,
that $\Ext^3(R_p/I^{[p^e]},R_p)=0$ for all $e$ by the flatness of
the Frobenius. This implies that $H^3_I(R_p)=0$ as well since in
positive characteristic $H^i_I(R_p) = \lim
\Ext^i(R_p/I^{[p^e]},R_p)$. On the other hand, it is well known that
in zero characteristic, $H^3_I(R_\QQ)$ is not zero. Hence this
provides an example where for the characteristic zero model we have
$\lambda_{0,3}=\lambda_{2,4} \neq 0$ whereas in all positive
characteristics $\lambda_{0,3}=\lambda_{2,4} = 0$.
\end{example}
The next example even shows that the vanishing of $\lambda_{a,i}$
can vary in an arithmetic progression:
\begin{example}
    Let $A$ be the homogeneous coordinate ring of $\PP^1 \times
    E$ where $E$ is the elliptic curve $\Proj
    \frac{\ZZ[x,y,z]}{x^3+y^3+z^3}$. Then $A$ is given as the quotient of
    $R$ (as above) by the ideal
    \[
        I=(\delta_1,\delta_2,\delta_3,x^3+y^3+z^3,ux^2+vy^2+zw^2,u^2x+v^2y+w^2z,u^3+v^3+w^3).
    \]
    The resolution of $A=R/I$ can be computed to be equal to
    \[
        0\to R \to R^6 \to R^{11} \to R^7 \to R \to 0
    \]
    and one verifies that $H^4_I(R_p) = 0$ if and only if $\op{char} k \equiv 2 \mod
    3$ (this essentially follows from the fact that depending on the
    modulus of $p$ mod $3$ the elliptic curve is supersingular or
    not). Hence we have that $\lambda_{0,2}=\lambda_{2,3} = 0$ if and only
    if $\op{char} k \equiv 2 \mod 3$.
\end{example}

%
%    Bibliographic information   %%%%%%%%%%%%%%%%%%%%%%%%%%%%%%%%%%%%%%
%

\providecommand{\bysame}{\leavevmode\hbox
to3em{\hrulefill}\thinspace}
\providecommand{\MR}{\relax\ifhmode\unskip\space\fi MR }
\providecommand{\MRhref}[2]{%
  \href{http://www.ams.org/mathscinet-getitem?mr=#1}{#2}
} \providecommand{\href}[2]{#2}

%
%    Index generation (if applicable)     %%%%%%%%%%%%%%%%%%%%%%%%%%%%%
%

\end{document}